\documentclass[11pt]{article}
\usepackage{graphicx}
\usepackage{amsmath}
\usepackage{amssymb}
\usepackage{amsthm}
\usepackage{appendix}
\usepackage{url}
\newcommand\RR{\mathbb{R}}

\newcommand\al\alpha
\newcommand\be\beta
\newcommand\de\delta
\newcommand\ep\varepsilon
\newcommand\tha\theta
\newcommand\ka\kappa
\newcommand\la\lambda
\newcommand\om\omega

\newcommand\iy\infty
\newcommand\pa\partial

\newcommand{\hyp}[5]{\,\mbox{}_{#1}F_{#2}\!\left(\genfrac{}{}{0pt}{}{#3}{#4};#5\right)}

\renewcommand\Re{\operatorname{Re}}

\numberwithin{equation}{section}
\newtheorem{theorem}{Theorem}[section]

\newtheorem{Remark}[theorem]{Remark}
\newenvironment{remark}{\begin{Remark}\rm}{\end{Remark}}

\begin{document}

\title{The two-dimensional fractional orthogonal derivative.}
\author{Enno Diekema \footnote{email adress: e.diekema@gmail.com}}
\maketitle

\begin{abstract}
\noindent
This paper is an edited and shortened version of Chapter 6 from the thesis of the author \cite{1}. First the one dimensional orthogonal derivative will be extended to the two-dimensional case. In the two-dimensional case we have to define the region of integration. In this paper we treat the integration over the square region and over the triangle region where in the last case we use biorthogonal polynomials expressed in terms of Appell  functions. Next the two-dimensional orthogonal derivative will be extended to the two-dimensional fractional orthogonal derivative. The results are highly dependent on the  $F_1$,\ $F_2$ and $F_3$ Appell functions.
\end{abstract}

\section{Introduction}
\setlength{\parindent}{0cm}
In Chapter 2 of the thesis of the author \cite[Ch. 2]{1} the orthogonal derivative of a function of one variable was defined. Just as this definition used orthogonal polynomials in one variable, we can try to approximate a partial derivative of a function in several variables by using orthogonal polynomials in several variables.
In this paper we will consider the two-variable case. Similarly, we will extend the definition in \cite[Ch. 3]{1} of the fractional orthogonal derivative to the two-variable case.
\\[3mm]
The contents of this paper are as follows. 
\\[2mm]
In section 2 we apply the method of \cite[Ch. 2]{1} to the two-variable case. We want to approximate $\frac{\partial^n f(x,y)}{\partial x^{n-k} \partial y^k}$. Only in the trivial case of a direct product of two orthogonality measures in one variable we can use fully orthogonal polynomials in two variables (again direct products) for this purpose. For other orthogonality measures on $\mathbb{R}^2$ we have to work with biorthogonal polynomials. The most obvious regions for the supports of these measures are the triangular region and the disk. We will only treat the triangular region with
orthogonal measure already considered by Appell and yielding biorthogonal polynomials expressed in terms of Appell hypergeometric functions.

In section 3 we apply the theory of \cite[Ch. 3]{1} to the two-variable case, by
which we obtain a two-dimensional fractional orthogonal derivative.

In section 4 we briefly treat the trivial case of the fractional Jacobi derivative with weight function on the square region.

In section 5 we treat the fractional biorthogonal derivative with weight function on the triangular region.

In the appendix we give an overview of the functions we need for the derivation of the formulas for the fractional orthogonal derivatives. We treat the $F_1$, $F_2$, $F_3$ Appell functions, the $H_2$ Horn function and the $F_P$, $F_Q$ and the $F_{PR}$ functions of Olsson.

\section{The two-dimensional orthogonal derivatives}
In \cite[Ch. 2]{1} we derived a formula for the orthogonal derivative
\[
f^{(n)}(x)=\dfrac{k_{n}n!}{h_{n}}\lim\limits_{\delta \downarrow 0}\dfrac{1}{\delta ^{n}}
\int\nolimits_{\RR}f(x+\delta u)p_{n}(u)d\mu(u)
\]
where the $p_n$ are orthogonal polynomials. This formula can be considered, for instance, for Jacobi, Laguerre or Hermite polynomials. We can also use discrete orthogonal polynomials, for example the Hahn polynomials.

In the two-variable case we will work with Taylor series, as we did in the one-dimensional case, but now there is more freedom. Write the Taylor expansion of a function $f$ around $(x,y)$ as
\begin{equation}
f(x+\delta u,y+\delta v)
=\delta^n\sum\limits_{m=0}^n\sum\limits_{l=0}^m\dfrac{\partial^m}{\partial x^{m-l}\partial y^{l}}f(x,y)
\dfrac{u^{m-l}v^{l}}{(m-l)!\,l!}+o(\delta^n)
\label{4.2.1}
\end{equation}
where $\delta$ is small. We want to approximate the partial derivative $\frac{\partial^n}{\partial x^{n-k} \partial y^k} f(x,y)$. For this purpose we look for a polynomial $p_{n,k}(u,v)$ of degree $n$ such that, when both sides of \eqref{4.2.1} are multiplied by $p_{n,k}(u,v)$ and next both sides are integrated over $\mathbb{R}^2$ with respect to a suitable orthogonality measure $d\mu(u,v)$, all terms in the double summation except for the term with $(m,l)=(n,k)$ become zero. So first we want
\[
\int\nolimits_{\RR}\int\nolimits_{\RR}p_{n,k}(u,v)u^{m-l}v^{l}d\mu(u,v)=0 \qquad \text{if} \ \ m<n
\]
This forces $p_{n,k}$ to be in the ($n+1$)-dimensional space ${\cal V}_n$ of orthogonal polynomials of degree $n$ with respect to the measure $\mu$. Next we require that
\begin{equation}
\int\nolimits_{\RR}\int\nolimits_{\RR}p_{n,k}(u,v)u^{n-l}v^ld\mu(u,v)=0
\qquad \text{if} \ \ k\neq l
\label{4.2.3}
\end{equation}
This can be equivalently written as
\[
\int_{\mathbb{R}^2} p_{n,k}(u,v) q_{m,l}(u,v) d\mu(u,v)=0 
\qquad \text{if} \ \ (m,l)\neq (n,k)
\]
where the polynomials $q_{n,l}$ are elements of ${\cal V}_n$ such that $q_{n,l}(u,v)=u^{n-l}v^l$ + polynomial of degree less than $n$. The polynomials $q_{n,l}$ form the so-called {\em monomial} or {\em monic} basis of ${\cal V}_n$. The above condition on $p_{n,k}\in {\cal V}_n$ now defines the polynomials $p_{n,k}$ uniquely, up to a constant factor, as the basis of ${\cal V}_n$ which is {\em biorthogonal} to the monomial basis. Then we get from \eqref{4.2.1} and the properties of $p_{n,k}$ that
\begin{equation}
\dfrac{\partial^n}{\partial x^{n-k}\partial y^{k}}f(x,y)=
(n-k)!\,k!\lim\limits_{\delta \downarrow 0}\dfrac{1}{\delta^n}
\dfrac{\int\nolimits_{\RR}\int\nolimits_{\RR}f(x+\delta u,y+\delta v)p_{n,k}(u,v)d\mu(u,v)}
{\int\nolimits_{\RR}\int\nolimits_{\RR}p_{n,k}(u,v)u^{n-k}v^k d\mu(u,v)}
\label{4.2.4}
\end{equation}
under some conditions on $f$, for instance that $f$ is a $C^n$ function on some neighbourhood of $(x,y)$ and, if $\mu$ has unbounded support, that $|f(u,v)|$ tends sufficiently fast to $0$ as $(u,v)$ tends to $\infty$. More generally, we call the right-hand side of \eqref{4.2.4} a {\em two-dimensional partial orthogonal derivative} of $f$ at $(x,y)$ if the limit exists. But we will in particular be interested in the numerator on the right-hand side (after the limit sign) for fixed $\delta$. We will only treat the square and the triangular region. For the disk see \cite[Section 5.2.2 with $d=2$]{6}, \cite[Chapter 3]{154}.

\subsection{The two-dimensional Jacobi derivative with weight function on the square region}
An easy example of \eqref{4.2.4} can be obtained by letting $d\mu(x,y)=d\mu_1(x) d\mu_2(y)$ and $p_{n,k}(x,y)=r_{n-k}(x)s_k(y)$ with the $r_n$ and $s_n$ one-variable orthogonal polynomials for orthogonality measures $\mu_1$ and $\mu_2$, respectively. Then the polynomials $p_{n,k}$ are two-variable orthogonal polynomials for the orthogonality measure $\mu$ while also $p_{n,k}(x,y)=$C$\,x^{n-k}y^k$ + polynomial of degree less than $n$. So the $p_{n,k}$ form the monomial basis, but this basis is also orthogonal, This is also biorthogonal to the monomial basis, by which the $p_{n,k}$ are the polynomials desired in \eqref{4.2.4}. Now the denominator after the limit sign on the right-hand side of \eqref{4.2.4} factorizes as the product
\[
\int_{\mathbb{R}} r_{n-k}(u) u^{n-k} d\mu_1(u) \times
\int_{\mathbb{R}} s_k(v) v^k d\mu_2(v)
=\frac{h_{n-k}'}{k_{n-k}'} \frac{h_{k}''}{k_{k}''}
\]
by \cite[(2.2.4)]{1}, where the accented and double accented $h_n$'s and $k_n$'s are in the obvious way, following \cite[(2.2.3)]{1}, related to the orthogonal polynomials $r_n$ and $s_n$, respectively. Specialization of \eqref{4.2.4} gives:
\begin{equation}
\dfrac{\partial ^{m+n}}{\partial ^{m}x\partial^{n}y}f(x,y)=
\dfrac{k_{m}\,m!}{h_{m}}\dfrac{k_{n}\,n!}{h_{n}}\lim\limits_{\delta	\rightarrow 0}\dfrac{1}{\delta ^{m+n}}\int\nolimits_{-1}^{1}\int\nolimits_{-1}^{1}f(x+\delta v,y+\delta u) p_{m}(u)q_{n}(v) d\mu(u)
d\mu(v)
\label{4.2.5}
\end{equation}
This is a very trivial case. If the measures $\mu_1$ and $\mu_2$ are supported on the segment $[-1,1]$ then the polynomials $p_{n,k}$ are orthogonal on the square $[-1,1]\times[-1,1]$. For example, take the product of the Jacobi measures $d\mu_1(x)=(1-x)^\alpha(1+x)^\beta dx$ and $d\mu_2(y)=(1-y)^\gamma(1+y)^\delta dy$. Substitution in \eqref{4.2.5} gives at last
\begin{multline*}
\dfrac{\partial ^{n+k}}{\partial x^{n}\partial y^{k}}f(x,y)=\lim\limits_{\delta	\rightarrow 0}\dfrac{1}{\delta^{n+k}}
\dfrac{k_{n,\alpha ,\beta }\ k_{k,\gamma,\delta }}{h_{n,\alpha ,\beta }\ h_{k,\gamma ,\delta }} \\
\times\int\nolimits_{-1}^{1}\int\nolimits_{-1}^{1}f(x+\delta u,y+\delta v) P_{n}^{(\alpha ,\beta) }
(u)P_{k}^{(\gamma,\delta)}(v)(1-u)^{\alpha }(1+u)^{\beta }(1-v) ^{\gamma }(1+v) ^{\delta}\,du\,dv
\end{multline*}

\subsection{The two-dimensional biorthogonal derivative with weight function on the triangular region}
The triangular region of ${\RR}^{2}$ is defined as $T^{2}:=\left\{(x,y) :0\leq x,y \le 1  \wedge x+y\le 1\right\}$,
on which we consider the Jacobi weight function defined as:
\[
W_{\alpha ,\beta ,\gamma }(x,y):=\dfrac{1}{B(\alpha +1,\beta +1,\gamma +1)}x^{\alpha }
y^{\beta}(1-x-y)^{\gamma }\qquad \alpha ,\beta ,\gamma >-1
\] 
which is normalized in such a way that its integral over $T^{2}$ is $1$. The function $B(x,y,z)$ is defined as an extension of the Beta function
\[
B(x,y,z):=\dfrac{\Gamma(x)\Gamma(y)\Gamma(z)}{\Gamma(x+y+z)}
\]
It is easy to prove
\[
B(x,y,z)=B(x,y+z)B(y,z)=B(y,x+z)B(x,z)=B(z,x+y)B(x,y)
\]

Let
\[
\left\langle f,g\right\rangle _{\alpha ,\beta ,\gamma}=\int\nolimits_{T^{2}}f(x,y) g(x,y)
W_{\alpha,\beta ,\gamma }(x,y)\;dx\;dy
\]
Classical orthogonal polynomials of one variable all satisfy a Rodrigues' formula. For obtaining something similar in connection with the triangular region Appell
considered the following analogue of Rodrigues' formula \cite[Section 2.4]{6}:
\begin{equation}
U_{k,n}^{\alpha ,\beta ,\gamma }(x,y) =
\big[x^\alpha y^\beta (1-x-y)^\gamma \big] ^{-1}\dfrac{\partial ^{n}}
{\partial x^{k}\partial y^{n-k}}\left[ x^{k+\alpha }y^{n-k+\beta }(1-x-y)^{n+\gamma } \right]
\qquad 0\leq k\leq n
\label{4.2.15}
\end{equation}
The set $\{U_{k,n}^{\alpha ,\beta ,\gamma }(x,y):0\leq k\leq n\}$ is a basis of
$\mathcal{V}_n^2 (W_{\alpha,\beta,\gamma})$ where $\mathcal{V}_n^2$\, is the space of orthogonal polynomials of degree $n$. See \cite[Prop. 2.4.3]{6}.

To write $U_{k,n}^{\alpha ,\beta ,\gamma}(x,y)$ as an $F_2$ Appell function we can use \cite[5.13(1)]{28}. There follows:
\begin{multline*}
U_{k,n}^{\alpha ,\beta ,\gamma }(x,y) =(\alpha+1)_k(\beta+1)_{n-k}(1-x-y)^n \\
F_2\left( -\gamma -n;-k,-n+k;\alpha +1,\beta +1;\dfrac{x}{x+y-1},\dfrac{y}{x+y-1}\right)
\end{multline*}
From this expression we see that $U_{k,n}^{\alpha ,\beta ,\gamma }(x,y)$ are polynomials of degree less or equal $n$. Using \cite[5.11(8)]{28} we get the much simpler form
\begin{multline}
U_{k,n}^{\alpha ,\beta ,\gamma }(x,y)
=(\alpha+1)_k(\beta+1)_{n-k}(1-x-y)^{-\gamma} \\
F_2(-\gamma -n;\alpha +1+k,\beta +1+n-k;\alpha +1,\beta+1;x,y)
\label{4.2.16}
\end{multline}
This basis is not orthogonal. It can be shown that $\left\langle U_{k,n}^{\alpha ,\beta,\gamma },U_{j,n}^{\alpha ,\beta ,\gamma }\right\rangle _{\alpha ,\beta,\gamma }\neq 0$ for $j\neq k$. If $n\neq m$ then $\left\langle U_{k,n}^{\alpha ,\beta,\gamma },U_{j,m}^{\alpha ,\beta,\gamma }\right\rangle _{\alpha ,\beta,\gamma }= 0$.

Let $V_{m,n}^{\alpha ,\beta ,\gamma }$ be defined with $0\leq m\leq n$ by
\begin{multline}
V_{m,n}^{\alpha ,\beta ,\gamma }(x,y):= \\
=\sum\limits_{i=0}^{m}\sum\limits_{j=0}^{n}(-1) ^{m+n+i+j}\dbinom{m}{i}\dbinom{n}{j}
\dfrac{(\alpha+1)_{m}(\beta +1)_n(\alpha+\beta+\gamma+2) _{m+n+i+j}}
{(\alpha+1)_{i}(\beta+1)_{j}(\alpha+\beta+\gamma+2)_{2m+2n}}x^{i}y^{j}
\end{multline}
This function can be written as an $F_2$ Appell function.
\begin{multline}
V_{m,n}^{\alpha ,\beta ,\gamma }(x,y) =(-1)^{n+m}
\dfrac{(\alpha+1)_{m}(\beta+1)_{n}(\alpha+\beta+\gamma+2) _{m+n}}
{(\alpha+\beta+\gamma+2) _{2n+2m}}\\
F_2(\alpha+\beta+\gamma+2+m+n;-m,-n;\alpha+1,\beta+1;x,y)
\label{4.2.17}
\end{multline}
Then $V_{m,n}$ is the orthogonal projection of $x^{m}y^{n}$ in $L^{2}\left\langle T^{2},W_{\mu }\right\rangle $ onto the space $\mathcal{V}_{m+n}$ of orthogonal polynomials of degree $m+n$ (so they are a monomial basis) and the two families of polynomials are biorthogonal with
\begin{equation}
\left\langle U_{k,n}^{\alpha ,\beta ,\gamma },V_{j,n-j}^{\alpha,\beta ,\gamma }\right\rangle _{\alpha ,\beta ,\gamma }=(-1)^n
\dfrac{(\alpha+1)_k(\beta+1)_{n-k}(\gamma+1)_n k!\,(n-k)!}
{(\alpha+\beta+\gamma+3)_{2n}}\delta _{j,k}\qquad 0\leq j,k\leq n
\label{4.2.17a}
\end{equation}
We can write this last inner product more generally by replacing the parameter $n-j$ by $m-j$. There follows
\[
\left\langle U_{k,n}^{\alpha ,\beta ,\gamma },V_{j,m-j}^{\alpha	,\beta ,\gamma }\right\rangle _{\alpha ,\beta ,\gamma }=0\qquad m\neq n
\]
The basis of polynomials $V_{m,n}^{\alpha,\beta,\gamma}$ as well as the above biorthogonality is already given by Appell for $\gamma=0$ and by Fackerel \& Littler  \cite{71} for the general case. If we compare with the notation for the general case at the beginning of Section 2 then we see that $U_{k,n}^{\alpha,\beta,\gamma}=p_{n,n-k}$ and $V_{l,m-l}^{\alpha,\beta,\gamma}=q_{m,m-l}$. Hence, for the partial orthogonal derivative \eqref{4.2.4} we have to use the $U$-basis in the present case.

\

Using the triangular region $T^2=\left\{(x,y) :0\leq x,y,x+y<1\right\}$ we get for the denominator of \eqref{4.2.4}
\begin{align}
I(m,l) \nonumber
&=\iint\nolimits_{0\leq u,v,u+v<1}U_{k,n}^{\alpha ,\beta ,\gamma }(u,v)
u^{m-l} v^l W_{\alpha,\beta,\gamma}(u,v) du dv= \\ \nonumber
&=\big[B(\alpha +1,\beta +1,\gamma +1)\big]^{-1} \\ 
&\qquad \qquad \times\int\nolimits_{0}^{1}v^{l}\int\nolimits_{0}^{1-v}u^{m-l}\dfrac{\partial ^{k}}
{\partial u^{k}}\dfrac{\partial^{n-k}}{\partial v^{n-k}}\left[u^{k+\alpha }v^{n-k+\beta }
(1-u-v)^{n+\gamma } \right]\,du\,dv
\label{4.2.4.3b}
\end{align} 
By integration by parts with respect to $u$ we see that $I\left( m,l\right)=0$ unless $l+k\leq m$. By integration by parts with respect to $v$ we see that $I(m,l) =0$ unless $n\leq l+k$. Together with $m\leq n$ there follows: $m=n$ and $l=n-k$. Substitution in \eqref{4.2.4.3b} gives
\begin{align}
I(n,n-k)
&B(\alpha +1,\beta +1,\gamma +1)= \nonumber \\
&=\int\nolimits_{0}^{1}v^{n-k}\int\nolimits_{0}^{1-v}u^{k}\dfrac{\partial ^{k}}
{\partial u^{k}}\dfrac{\partial ^{n-k}}{\partial v^{n-k}}\left[ u^{k+\alpha }v^{n-k+\beta }(1-u-v)^{n+\gamma } \right]\,du\,dv
\label{4.2.4.3c}
\end{align}
Let
\[
I_{1}=\int\nolimits_{0}^{1-v}u^{k}\dfrac{\partial ^{k}}{\partial u^{k}}
\left[\dfrac{\partial ^{n-k}}{\partial v^{n-k}}\left[ u^{k+\alpha }v^{n-k+\beta }(1-u-v)^{n+\gamma } \right]\right]\,du
\]
After $k$-fold integration by parts there remains
\[
I_{1}=(-1)^k\Gamma(k+1)\int\nolimits_{0}^{1-v}u^{k+\alpha}\dfrac{\partial ^{n-k}}{\partial v^{n-k}}\left[ v^{n-k+\beta }(1-u-v) ^{n+\gamma }\right] du
\]
Substitution in \eqref{4.2.4.3c} gives
\begin{align}
I(n,n-k)
&\big[(-1)^k\Gamma(k+1)\big]^{-1} B(\alpha +1,\beta +1,\gamma +1)= \nonumber \\
&=\int\nolimits_{0}^{1}v^{n-k}\int\nolimits_{0}^{1-v}u^{k+\alpha}
\dfrac{\partial ^{n-k}}{\partial v^{n-k}}\left[v^{n-k+\beta }(1-u-v)^{n+\gamma } \right]\,du\,dv \nonumber \\
&=\int\nolimits_{0}^{1}u^{k+\alpha}\int\nolimits_{0}^{1-u}v^{n-k}
\dfrac{\partial ^{n-k}}{\partial v^{n-k}}\left[v^{n-k+\beta }(1-u-v)^{n+\gamma } \right]\,dv\,du
\label{4.2.4.7}
\end{align}
Let
\[
I_{2}=\int\nolimits_{0}^{1-u}v^{n-k}
\dfrac{\partial ^{n-k}}{\partial v^{n-k}}\left[v^{n-k+\beta }(1-u-v)^{n+\gamma } \right]\,dv
\] 
After $(n-k)$-fold integration by parts we get
\[
I_{2}=(-1)^{n-k}\Gamma(n-k+1)\int\nolimits_{0}^{1-u}v^{n-k+\beta}(1-u-v)^{n+\gamma}dv
\] 
The integral is a Beta function.
\[
\int\nolimits_{0}^{1-u}v^{n-k+\beta}(1-u-v)^{n+\gamma}dv=
(1-u)^{k+2n+\beta +\gamma +1}B(n-k+\beta +1,n+\gamma +1)
\]
with convergence conditions: $\beta,\gamma>-1$. So for $I_2$ we get
\[
I_{2}=(-1)^{n-k}\Gamma(n-k+1)(1-u)^{k+2n+\beta +\gamma +1}B(n-k+\beta +1,n+\gamma +1)
\]
Substitution in \eqref{4.2.4.7} gives
\begin{multline*}
I(n,n-k)\big[(-1)^n\Gamma(k+1)\Gamma(n-k+1)\big]^{-1} B(\alpha +1,\beta +1,\gamma +1)= \\
=B(n-k+\beta +1,n+\gamma +1)
\int\nolimits_{0}^{1}u^{k+\alpha}(1-u)^{k+2n+\beta +\gamma +1}du
\end{multline*}
The integral is again a Beta function. So we get
\begin{align*}
I(n,n-k)\big[(-1)^n
&\Gamma(k+1)\Gamma(n-k+1)\big]^{-1}B(\alpha +1,\beta +1,\gamma +1)= \\
&=B(n-k+\beta +1,n+\gamma +1)B(k+\alpha+1,2n-k+\beta+\gamma+2) \\
&=B(k+\alpha +1,n-k+\beta+1,n+\gamma +1)
\end{align*}
with convergence conditions: $\Re{\alpha},\Re{\beta},\Re{\gamma}>-1$. After some  manipulations we get
\[
I(n,n-k)=(-1)^n\dfrac{(\alpha+1)_k(\beta+1)_{n-k}(\gamma+1)_n k!\,(n-k)!}
{(\alpha+\beta+\gamma+3)_{2n}}
\]
and this is equal to the constant in \eqref{4.2.17a}. This also follows from the fact that
\begin{align*}
I(n,n-k)
&=<U_{k,n},\,u^k\,v^{n-k}> \\
&=<U_{k,n},\,V_{k,n-k}+\text{polynomial of degree less than $n$}> \\
&=<U_{k,n},\,V_{k,n-k}>+<U_{k,n},\text{polynomial of degree less than $n$}> \\
&=<U_{k,n},\,V_{k,n-k}>
\end{align*}
Substitution in \eqref{4.2.4} gives
\begin{multline}
\dfrac{\partial ^{n}}{\partial x^{k}\partial y^{n-k}}f(x,y) =
\dfrac{(-1)^n}{B(\alpha+1+k,\beta+1+n-k,\gamma+1+n)} \\
\lim\limits_{\delta \downarrow 0}\dfrac{1}{\delta ^{n}}\iint\nolimits_{0\leq u,v,u+v<1}
f(x+\delta u,y+\delta v)
U_{k,n}^{\alpha ,\beta ,\gamma }(u,v) u^{\alpha }v^{\beta}(1-u-v) ^{\gamma }\,dv\,du
\label{4.6.14a}
\end{multline}
with convergence conditions: $\alpha,\beta,\gamma>-1$. Using \eqref{4.2.16} gives
\begin{multline*}
\dfrac{\partial ^{n}}{\partial x^{k}\partial y^{n-k}}f(x,y) =
\dfrac{(-1)^n\Gamma(\alpha+\beta+\gamma+3+2n)}{\Gamma(\alpha+1)\Gamma(\beta+1)
\Gamma(\gamma+1+n)} \\
\lim\limits_{\delta \downarrow 0}\dfrac{1}{\delta ^{n}} 
\iint\nolimits_{0\leq u,v,u+v<1}f(x+\delta u,y+\delta v)
F_2\left( 
\begin{array}{c}
-\gamma -n;\alpha +1+k,\beta +1+n-k \\ 
\alpha +1,\beta+1%
\end{array}%
;
u,v\right) 
u^{\alpha }v^{\beta}\,dv\,du
\end{multline*}

\section{The two-dimensional fractional orthogonal derivative}
We have seen that for the orthogonal derivative in two dimensions there are several choices of region and orthogonality measures for which explicit results can be obtained. For the fractional derivative in two dimensions we have the same possibilities. So we treat first the trivial case of the fractional orthogonal derivative with weight function on the square. Then we treat the fractional biorthogonal derivative with weight function on the triangular region.

Analogous to the Weyl fractional integral in one dimension we define the Weyl fractional integral in two dimensions as:
\[
W^{-\mu ,-\nu }[f](x,y):=\dfrac{1}{\Gamma(\mu)\Gamma(\nu) }
\int\nolimits_{x}^{\infty}\int\nolimits_{y}^{\infty }f(U,V)(U-x)^{\mu -1}(V-y)^{\nu -1}dVdU
\]
The fractional partial derivative will be given by
\begin{multline}
W^{\mu,\nu}[f](x,y)=(-1)^m\dfrac{\partial^m}{\partial x^{m-l}\partial y^l}
W^{\mu -m+l,}{}^{\nu -l}[f](x,y) \\
=\dfrac{(-1)^m}{\Gamma(m-l-\mu)\Gamma(l-\nu)}
\dfrac{\partial ^m}{\partial x^{m-l}\partial y^l}
\left[\int\nolimits_{x}^{\infty}\int\nolimits_{y}^{\infty }f(U,V)(U-x)^{m-l-\mu -1}
(V-y)^{l-\nu -1}dVdU\right]
\label{4.3.1.4}
\end{multline}
Analogous to \cite[(3.3.1)]{1} and in view of \eqref{4.2.4} we can write
\[
W^{\mu ,\nu }[f](x,y) =\lim\limits_{\delta \downarrow 0}W^{\mu,\nu,m,l}_\delta [f](x,y)
\]
where
\begin{align}
W^{\mu,\nu,m,l}_\delta [f](x,y):&=(-1)^mD^{m,l}_\delta \big[ W^{\mu-m+l,\nu-l}[f]\big](x,y) \nonumber \\
&=(-1)^mW^{\mu-m+l,\nu-l}\big[ D^{m,l}_\delta [f]\big](x,y)
\label{6.4.3}
\end{align}
with
\begin{equation}
D^{m,l}_\delta [g](x,y):=\dfrac{(m-l)!\,l!}{\delta^m}\dfrac{}{}\dfrac{\int\nolimits_{\RR}\int\nolimits_{\RR}g(x+\delta u,y+\delta v)p_{m,l}(u,v)d\mu(u,v)}
{\int\nolimits_{\RR}\int\nolimits_{\RR}p_{m,l}(u,v)u^{m-l}v^l d\mu(u,v)}
\label{6.4.4}
\end{equation}
After application of \eqref{6.4.3} and \eqref{6.4.4} to \eqref{4.3.1.4} we obtain
\begin{multline*}
W^{\mu,\nu,m,l}_\delta [f](x,y)=\dfrac{(-1)^m}{\Gamma(m-l-\mu)\Gamma(l-\nu)}
\dfrac{(m-l)!\,l!}{\int\nolimits_{\RR}\int\nolimits_{\RR}p_{m,l}(u,v)u^{m-l}v^l d\mu(u,v)}
\dfrac{1}{\delta^m} \\
\times\int\nolimits_{\RR}\int\nolimits_{\RR}\left[\int\nolimits_{x+\delta u}^{\infty}
\int\nolimits_{y+\delta v}^{\infty }f(U,V)(U-x-\delta u)^{m-l-\mu -1}(V-y-\delta v)^{l-\nu -1}dVdU \right] \\
\times p_{m,l}(u,v) d\mu(u,v)
\end{multline*}
With $U=x+\delta s$ and $V=y+\delta t$ there follows
\begin{multline*}
W^{\mu,\nu,m,l}_\delta [f](x,y)=\dfrac{(-1)^m}{\Gamma(m-l-\mu)\Gamma(l-\nu)}
\dfrac{(m-l)!\,l!}{\int\nolimits_{\RR}\int\nolimits_{\RR}p_{m,l}(u,v)u^{m-l}v^l d\mu(u,v)}
\dfrac{1}{\delta^{\mu+\nu}} \\
\times\int\nolimits_{\RR}\int\nolimits_{\RR}\left[\int\nolimits_{x}^{\infty}\int\nolimits_{y}^{\infty }f(x+\delta s,y+\delta t)(s-u)^{m-l-\mu -1}(t-v)^{l-\nu -1}dtds \right]
p_{m,l}(u,v) d\mu(u,v)
\end{multline*}
Interchange of the  inner double integral with the outer double integral gives
\begin{multline}
W^{\mu,\nu,m,l}_\delta [f](x,y)=\dfrac{(-1)^m}{\Gamma(m-l-\mu)\Gamma(l-\nu)}
\dfrac{(m-l)!\,l!}
{\int\nolimits_{\RR}\int\nolimits_{\RR}p_{m,l}(u,v)u^{m-l}v^l d\mu(u,v)}
\dfrac{1}{\delta^{\mu+\nu}} \\
\times\int\nolimits_{\RR}\int\nolimits_{\RR}f(x+\delta s,y+\delta t)
\left[\int\nolimits_{u=-\infty}^s\int\nolimits_{v=-\infty}^t (s-u)^{m-l-\mu-1}(t-v)^{l-\nu-1} p_{m,l} (u,v)d\mu(u,v)\right]dsdt
\label{6.4.5}
\end{multline}
In subsequent subsections we will write \eqref{6.4.5} in more explicit form for the cases of the square and the triangular regions.

\section{The two-dimensional fractional Jacobi derivative with weight function on the square.}
Let $d\mu(x,y)=d\mu_1(x)d\mu_2(y)$ and $p_{n,k}(x,y)=r_{n-k}(x)s_k(y)$ as in the beginning of section 2.1. Then by \eqref{4.2.5} formula \eqref{6.4.5} takes the form
\begin{multline*}
W^{\mu,\nu,m,l}_\delta [f](x,y)=\dfrac{(-1)^m(m-l)!\ l!\ k'_{m-l}\ k''_l}{\Gamma(m-l-\mu)\Gamma(l-\nu)h'_{m-l}\ h''_l}
\dfrac{1}{\delta^{\mu+\nu}}
\int\nolimits_{\RR}\int\nolimits_{\RR}f(x+\delta s,y+\delta t) \\
\times\left[\int\nolimits_{u=-\infty}^s(s-u)^{m-l-\mu-1}r_{m-l}(u)d\mu_1(u)\right]
\left[\int\nolimits_{v=-\infty}^t(t-v)^{l-\nu-1}s_l(v)d\mu_2(v)\right]dsdt
\end{multline*}
Now specialize to $d\mu_1(x)=(1-x)^\alpha(1+x)^\beta dx$ on $[-1,1]$ for $r_n(x)=P^{(\alpha,\beta)}_n(x)$ and $d\mu_2(y)=(1-y)^\gamma(1+y)^\delta\ dy$ for $s_n(y)=P^{(\gamma,\delta)}_n(y)$. Then we get for the fractional Jacobi derivative
\begin{multline}
W^{\mu,\nu,m,l}_\delta [f](x,y)
=\dfrac{(-1)^m(m-l)!\ l!\ k'_{m-l}\ k''_l}{\delta^{\mu+\nu}\  h'_{m-l}\ h''_l} \\
\times \int\nolimits_{-1}^1\int\nolimits_{-1}^1f(x+\delta s,y+\delta t)
\left[\int\nolimits_{-\infty}^s(s-u)^{m-l-\mu-1}P^{(\alpha,\beta)}_{m-l}(u)(1-u)^\alpha(1+u)^\beta du\right] \\
\times \left[\int\nolimits_{-\infty}^t(t-v)^{l-\nu-1}P^{(\gamma,\delta)}_l(v)(1-v)^\gamma(1+v)^\delta\ dv\right] dsdt \qquad
\label{4.3.1.10}
\end{multline}
We now look at the two-dimensional analogue of \cite[(3.3.6)]{1}. There we had transformed the 2-fold integral into a sum of two single integrals, where the integrands are a product of the function $f(x+\delta s)$ with a factor depending on a convolution integral of an orthogonal polynomial. In the two-dimensional case we can transform the $4$-fold integral in such a manner that there results a sum of four $2$-fold integrals with a product of the function $f(x+\delta s,y+\delta t)$ and a factor depending on the $2$-fold convolution integral of the orthogonal polynomials. In the special case of a product of Jacobi polynomials we can use the same computation as in the one-dimensional case. We give here only the result.
\begin{align*} 
W^{\mu ,\nu ,m+n,n}[f](x,y)
&\left[(-1)^{m+n}\dfrac{k_{m}m!}{h_{m}}\dfrac{k_{n}n!}{h_{n}}\dfrac{1}{\delta^{\mu +\nu }}
\dfrac{1}{\Gamma(m-\mu)\Gamma(n-\nu) }\right] ^{-1}=  \\
&=\int\nolimits_{-1}^{1}\int\nolimits_{-1}^{1}f(x+\delta s,y+\delta t)J_{1}(t,n,\alpha ,\beta ,\nu )
J_{1}(s,m,\gamma ,\delta ,\mu) \,dt\,ds+  \\
&+\int\nolimits_{1}^{\infty}\int\nolimits_{-1}^{1}f(x+\delta s,y+\delta t)J_{1}(t,n,\alpha,\beta,\nu) J_{2}(s,m,\gamma ,\delta ,\mu)\,dt\,ds+ \\
&+\int\nolimits_{-1}^{1}\int\nolimits_{1}^{\infty}f(x+\delta s,y+\delta t)J_{2}(t,n,\alpha,\beta,\nu) J_{1}(s,m,\gamma ,\delta ,\mu) \,dt\,ds+ \nonumber \\
&+\int\nolimits_{1}^{\infty}\int\nolimits_{1}^{\infty}f(x+\delta s,y+\delta t)
J_{2}(t,n,\alpha,\beta,\nu)J_{2}(s,m,\gamma ,\delta,\mu) \,dt\,ds
\end{align*}
The integrals $J_{1}$ and $J_{2}$ are known \cite[(3.5.9)]{1} and \cite[(3.5.10)]{1}. We obtain
\begin{multline*}
J_{1}(\xi ,\lambda ,\alpha ,\beta ,\nu) =(-1) ^{\lambda }\dfrac{\Gamma(\lambda +\beta +1)
	\Gamma(\lambda-\nu)}{2^{\lambda-\nu }\lambda !\,\Gamma(\lambda-\nu+\beta+1)}
(1-\xi)^{\lambda+\alpha-\nu}(1+\xi)^{\lambda +\beta -\nu } \\
\hyp21{-\nu ,2\lambda -\nu +\alpha +\beta +1}{\lambda -\nu+\beta +1}{\dfrac{1+\xi }{2}}
\end{multline*}
\begin{multline*}
J_{2}(\xi ,\lambda ,\alpha ,\beta ,\nu) =(-1) ^{\lambda }\dfrac{2^{\lambda +\alpha +\beta +1}}
{(\xi+1)^{\nu +1}}\dfrac{\Gamma(\lambda -\nu)}{\Gamma(-\nu) \lambda !}
\dfrac{\Gamma(\lambda +\alpha +1)\Gamma(\lambda +\beta +1)}{\Gamma(2\lambda +\alpha+\beta +2) } \\
\hyp21{\nu +1,\lambda +\beta +1}{2\lambda +\alpha +\beta +2}{\dfrac{2}{\xi +1}}
\end{multline*}

\section{The fractional biorthogonal derivative with weight function on the triangular region}
\begin{figure}[ht]
	\centering
	\includegraphics[height=4cm]{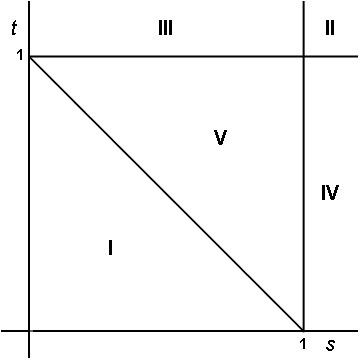}
	\caption{The five possible regions of the position of the point $(s,t)$.}
	\label{Figure 6.1a}
\end{figure}

\begin{figure}[ht]
	\includegraphics[height=3.5cm,width=16.8cm]{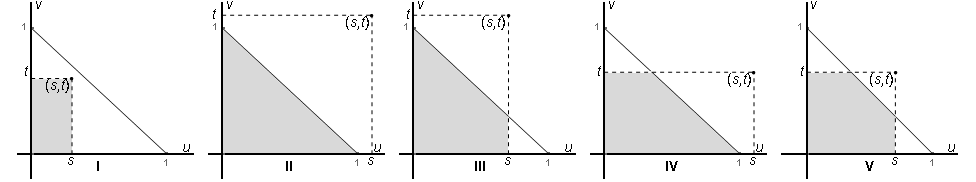}
	\caption{The five possibilities of the position of the point $(s,t)$ with respect to the $u,v$ axes.}
	\label{Figure 6.1}
\end{figure}

To develop a formula for the fractional biorthogonal derivative on the triangular region, we use \eqref{6.4.5} for the triangular region $T^2$.
\begin{multline}
W^{\mu,\nu,n,n-k}_\delta [f](x,y)\left[ \dfrac{h_{1}(\alpha ,\beta ,\gamma,k,n)}
{\Gamma(k-\mu)\Gamma(n-k-\nu)}\dfrac{1}{\delta ^{n}}\right] ^{-1}=
\int\nolimits_{\RR}\int\nolimits_{\RR}f(x+\delta s,y+\delta t) \\
\times\left[\iint\nolimits_{T^2 \cap \big( (-\infty,s]\times(-\infty,t] \big)} (s-u)^{k-\mu-1}(t-v)^{n-k-\nu-1} U^{\alpha,\beta,\gamma}_{k,n}(u,v)u^\alpha v^\beta (1-u-v)^\gamma dudv)\right]dsdt
\label{4.9.4a}
\end{multline}
with
\begin{equation}
h_{1}(\alpha ,\beta ,\gamma ,k,n)=\dfrac{1}{B(\alpha+1+k,\beta+1+n-k,\gamma+1+n)}
\label{4.9.4}
\end{equation}
We see that the integration domain $T^2\cap((-\infty,s]\times(-\infty,t])$ of the inner double integral in \eqref{4.9.4a} depends on the position of the point $(s,t)$ with respect to the triangle $\left\{ \left(u,v\right) \text{ }|\text{ }0<v<u<1\right\}$. Clearly $T^2\cap((-\infty,s]\times(-\infty,t])$ is empty if $(s,t)$ is outside the
first quadrant, while the first quadrant in the $(s,t)$ plane can be split up
as the union of five regions with disjoint interiors as in Figure \ref{Figure 6.1a} such that the corresponding regions $T^2\cap((-\infty,s]\times(-\infty,t])$ inside $T^2$ in the
$(u,v)$ plane are given by Figure \ref{Figure 6.1}.

$\text{I.     }$ The point $(s,t)$ lies inside the triangle with $0<s<1,0<t<1-s$. The appropriate part of the integral for $W$ becomes
\[
I=\int_{s=0}^{1}\int_{t=0}^{1-s}\left(\int_{v=0}^{t}\int_{u=0}^{s}(.)\,du\,dv\right) \,dt\,ds
\]
$\text{II.   }$ The point $(s,t)$ lies outside the triangle with
$1<s<\infty ,1<t<\infty$. The appropriate part of the integral for $W$ becomes
\[
I=\int_{s=1}^{\infty }\int_{t=1}^{\infty }\left(\int_{v=0}^{1}\int_{u=0}^{1-v}(.)\,du\,dv\right) \,dt\,ds
\]
$\text{III.    }$ The point $(s,t)$ lies outside the triangle with
$0<s<1<t<\infty$. The appropriate part of the integral for $W$ becomes
\[
I=\int_{t=1}^{\infty }\int_{s=0}^{1}\left(\int_{u=0}^{s}\int_{v=0}^{1-u}(.) \,dv\,du\right) \,ds\,dt
\]
$\text{IV.    }$ The point $(s,t)$ lies outside the triangle with $0<t<1<s<\infty$. The appropriate part of the integral for $W$ becomes
\[
I=\int_{s=1}^{\infty }\int_{t=0}^{1}\left(\int_{v=0}^{t}\int_{u=0}^{1-v}(.)\,du\,dv\right) \,dt\,ds
\]
$\text{V.    }$ The point $(s,t)$ lies outside the triangle with
$0<s<1,\ 1-s<t<1$. We divide the region of integration in two parts. The appropriate part of the integral for $W$ becomes
\[
I=\int_{s=0}^{1}\int_{t=1-s}^{1}\left(\int_{v=0}^{1-s}\int_{u=0}^{s}(.)\,du\,dv\right)\,dt\,ds
+\int_{s=0}^{1}\int_{t=1-s}^{1}\left( \int_{v=1-s}^{t}\int_{u=0}^{1-v}(.)\,du\,dv\right) \,dt\,ds
\]

For the fractional derivative we get (the roman numbers refer to the cases in Figure \ref{Figure 6.1})
\begin{align}
W^{\mu,\nu}
&[f](x,y)\left[ \dfrac{	h_{1}(\alpha ,\beta ,\gamma ,k,n) }{\Gamma(k-\mu)\Gamma(n-k-\nu) }
\dfrac{1}{\delta ^{\nu +\mu }}\right] ^{-1}= \nonumber \\
&=\int_{s=0}^{1}\int_{t=0}^{1-s}f(x+\delta s,y+\delta t)
\left( \int_{v=0}^{t}\int_{u=0}^{s}\text{ \ }(.)\,du\,dv\right) \,dt\,ds+\qquad \qquad \quad\  \text{(I)} \nonumber \\
&+\int_{s=1}^{\infty }\int_{t=1}^{\infty }f(x+\delta s,y+\delta t)
\left(\int_{v=0}^{1}\int_{u=0}^{1-v}(.)\,du\,dv\right) \,dt\,ds+\qquad \qquad \qquad\  \text{(II)} \nonumber \\
&+\int_{t=1}^{\infty }\int_{s=0}^{1}\text{  }f(x+\delta s,y+\delta t) \left(\int_{u=0}^{s}\int_{v=0}^{1-u}(.)\,dv\,du\right) \,ds\,dt+\qquad \qquad \quad\ \  \text{(III)} \nonumber \\
&+\int_{s=1}^{\infty }\int_{t=0}^{1}f(x+\delta s,y+\delta t)
\left(\int_{v=0}^{t}\int_{u=0}^{1-v}(.)\,du\,dv\right) \,dt\,ds+\qquad \qquad \qquad \text{   (IV)} \nonumber \\
&+\int_{s=0}^{1}\int_{t=1-s}^{1}f(x+\delta s,y+\delta t)
\left( \int_{v=0}^{1-s}\int_{u=0}^{s}(.)\,du\,dv\right) \,dt\,ds+\qquad \qquad\quad\  \text{ (V)} \nonumber \\ \nonumber
&\quad\quad\quad +\int_{s=0}^{1}\int_{t=1-s}^{1}f(x+\delta s,y+\delta t) \left(\int_{v=1-s}^{t}\int_{u=0}^{1-v}(.)\,du\,dv\right) \,dt\,ds \\
\label{4.9.11}
\end{align}
The abbreviation $(.)$ stands for 
\begin{equation}
(.)=(s-u)^{k-\mu -1}(t-v)^{n-k-\nu-1}U_{k,n}^{\alpha ,\beta ,\gamma }(u,v) u^{\alpha}v^{\beta }(1-u-v) ^{\gamma }
\label{4.9.12b}
\end{equation}
For $U_{k,n}^{\alpha ,\beta ,\gamma}(u,v) $ we repeat formulas \eqref{4.2.15} and \eqref{4.2.16}:
\begin{equation}
U_{k,n}^{\alpha ,\beta ,\gamma }(u,v) =\dfrac{1}{u^{\alpha }v^{\beta }(1-u-v) ^{\gamma }} \\
\dfrac{\partial ^{n}}{\partial u^{k}\partial v^{n-k}}\left[ u^{k+\alpha }v^{n-k+\beta }
(1-u-v)^{n+\gamma } \right] 
\label{4.9.12a}
\end{equation}
\begin{multline*}
\qquad\qquad \ U_{k,n}^{\alpha ,\beta ,\gamma }(x,y)
=(\alpha+1)_k(\beta+1)_{n-k}(1-x-y)^{-\gamma} \\
F_2(-\gamma -n;\alpha +1+k,\beta +1+n-k;\alpha +1,\beta+1;x,y)
\end{multline*}
for $0\leq k\leq n$. We can use both formulas to compute the fractional derivatives. The results will be of course the same, but with \eqref{4.9.12a} there appear some nice results and so we continue with \eqref{4.9.12a}.

\

For the integrals inside the brackets in \eqref{4.9.11} we can use \eqref{4.9.12b} and \eqref{4.9.12a} and write in general
\begin{multline*}
\iint_{V}(.)\,du\,dv= \\
=\iint_{V}\dfrac{\partial ^{k}}{\partial u^{k}}\dfrac{\partial ^{n-k}}
{\partial v^{n-k}}\left[ u^{k+\alpha}v^{n-k+\beta }(1-u-v)^{n+\gamma } \right](s-u)^{k-\mu -1}
(t-v)^{n-k-\nu-1}\,du\,dv
\end{multline*}
After $k$-fold integration by parts with respect to the variable $u$ we get
\begin{multline*}
\iint_{V}(.)\,du\,dv\left[(-1)^{k}\dfrac{\Gamma(k-\mu)}{\Gamma(-\mu)}\right]^{-1}= \\
=\iint\nolimits_{V}(s-u) ^{-\mu -1}\dfrac{\partial ^{n-k}}{\partial v^{n-k}}\left[ u^{k+\alpha }v^{n-k+\beta }(1-u-v)^{n+\gamma } \right](t-v)^{n-k-\nu -1}\,du\,dv
\end{multline*}
After $(n-k) $-fold integration by parts with respect to the variable $v$ we get
\begin{multline*}
\iint_{V}(.) \,du\,dv\left[ (-1) ^{n}\dfrac{\Gamma(k-\mu)}{\Gamma(-\mu) }
\dfrac{\Gamma(n-k-\nu)}{\Gamma(-\nu)}\right] ^{-1}= \\
=\iint\nolimits_{V}u^{k+\alpha }(s-u) ^{-\mu -1}(1-u-v) ^{n+\gamma }v^{n-k+\beta }
(t-v) ^{-\nu -1}\,du\,dv
\end{multline*}
Combination of \eqref{4.9.4} and \eqref{4.9.11} gives
\begin{align*}
W^{\mu,\nu}[f](x,y)
&\left[ \dfrac{1}{\Gamma(-\mu) \Gamma(-\nu)}
\dfrac{1}{B(k+\alpha+1,n+\gamma+1,n-k+\beta +1)}\dfrac{1}{\delta ^{\nu +\mu }}\right] ^{-1}= \\
&=\int_{0}^{1}\int_{0}^{1-s} f(x+\delta s,y+\delta t) I_{1}(s,t) \,dt\,ds+\qquad \qquad  \text{ \ \ \ \ (I)} \nonumber \\
&+\int_{1}^{\infty }\int_{1}^{\infty } f(x+\delta s,y+\delta t) I_{3}(s,t) \,dt\,ds+\qquad \qquad\qquad \text{(II)}  \\
&+\int_{1}^{\infty }\int_{0}^{1}\text{  } f(x+\delta s,y+\delta t) I_{2}(s,t) \,ds\,dt+\qquad \qquad \qquad \text{(III)} \\
&+\int_{1}^{\infty }\int_{0}^{1}\text{  } f(x+\delta s,y+\delta t) I_{4}(s,t) \,dt\,ds+\qquad \qquad \qquad \text{(IV)} \\
&+\int_{0}^{1}\int_{1-s}^{1}\text{  } f(x+\delta s,y+\delta t) I_{5}(s,t) \,dt\,ds+ \qquad \qquad \qquad \text{(V)} \\ 
&\qquad \qquad \qquad +\int_{0}^{1}\int_{1-s}^{1} f(x+\delta s,y+\delta t) I_{6}(s,t)\,dt\,ds 
\label{4.9.16}
\end{align*}

Putting $a=k+\alpha +1$, $b=n-k+\beta+1$, $c=-\mu$, $d=-\nu$ and $e=-n-\gamma$ we get for the functions $I_{1..6}(s,t)$
\begin{equation}
I_{1}(s,t)=\int_{0}^{t}\left(\int_{0}^{s}u^{a-1}(s-u)^{c-1}(1-u-v)^{-e}du\right)v^{b-1}(t-v) ^{d-1}dv
\label{4.9.17a}
\end{equation}
with region for $s$ and $t$:$\qquad 0<s<1,0<t<1-s$
\begin{equation}
I_{2}(s,t)=\int_{0}^{1}\left(\int_{0}^{1-v}u^{a-1}(s-u)^{c-1}(1-u-v)^{-e}du\right) v^{b-1}(t-v)^{d-1}dv
\label{4.9.17c}
\end{equation}
with region for $s$ and $t$:$\qquad 1<s<\infty ,1<t<\infty$
\begin{equation}
I_{3}(s,t) =\int_{0}^{s}\left( \int_{0}^{1-u}v^{b-1}(t-v)^{d-1}(1-u-v)^{-e}dv \right) u^{a-1}(s-u)^{c-1}du
\label{4.9.17b}
\end{equation}
with region for $s$ and $t$:$\qquad 0<s<1<t<\infty$
\begin{equation}
I_{4}(s,t)=\int_{0}^{t}\left(\int_{0}^{1-v}u^{a-1}(s-u)^{c-1}(1-u-v)^{-e}du\right) v^{b-1}(t-v)^{d-1}dv
\label{4.9.17d}
\end{equation}
with region for $s$ and $t$:$\qquad 0<t<1<s<\infty$
\begin{equation}
I_{5}(s,t)=\int_{0}^{1-s}\left(\int_{0}^{s}u^{a-1}(s-u)^{c-1}(1-u-v)^{-e}du\right) v^{b-1}(t-v)^{d-1}dv
\label{4.9.17e}
\end{equation}
\begin{equation}
I_{6}(s,t)=\int_{1-s}^{t}\left(\int_{0}^{1-v}u^{a-1}(s-u)^{c-1}(1-u-v)^{-e}du\right) v^{b-1}(t-v) ^{d-1}dv
\label{4.9.17f}
\end{equation}
with region for $s$ and $t$:$\qquad 0<s<1,\ 1-s<t<1$

We will treat the evaluation of the integrals $I_{1,\dots, 6}$ in the next subsections.

\subsection{The integral $I_1(s,t)$}
For the computation of the integral $I_{1}(s,t)$ we write the integral \eqref{4.9.17a} as
\begin{equation}
I_{1}(s,t)=s^{a+c-1}t^{b+d-1}\int_{0}^{1}\int_{0}^{1}u^{a-1}v^{b-1}(1-u)^{c-1}(1-v)^{d-1}( 1-su-tv) ^{-e}\,du\,dv
\label{4.9.19}
\end{equation}
The integrals are convergent if $0<\Re(a)<\Re(a+c)$ and $0<\Re(b)<\Re(b+d)$. Formula \eqref{4.9.19} is the integral representation for the $F_2$ Appell function \cite[5.8.1.(2)]{28}. We get
\[
I_{1}(s,t) =B(a,c)B(b,d)s^{a+c-1}t^{b+d-1}F_2( e;a,b;a+c,b+d;s,t)
\]
\subsection{The integral $I_2(s,t)$}
For the computation of the integral $I_{2}(s,t)$ we write the integral \eqref{4.9.17c} as
\begin{equation}
I_{2}(s,t)=s^{c-1}t^{d-1}\int_{0}^{1}\int_{0}^{1-v}u^{a-1}v^{b-1}(1-u-v)
^{-e}\left( 1-\dfrac{1}{s}u\right) ^{c-1}\left( 1-\dfrac{1}{t}v\right)^{d-1}\,du\,dv
\label{4.9.23}
\end{equation}
The integrals are convergent if $\Re(a),\Re(b),\Re(1-e)>0$. Formula \eqref{4.9.23} is the integral representation for the $F_3$ Appell function \cite[5.8.1.(3)]{28}. We get 
\begin{equation}
I_{2}(s,t)=B(a,b,1-e)s^{c-1}t^{d-1}
F_3\left( a,b;1-c,1-d;a+b+1-e;\dfrac{1}{s},\dfrac{1}{t}\right)
\label{4.9.24}
\end{equation}
\subsection{The integral $I_3(s,t)$}
For the computation of the integral $I_{3}(s,t)$ we write the integral \eqref{4.9.17b} as
\begin{multline}
I_{3}(s,t)=s^{a+c-1}t^{d-1} \\
\int_{0}^{1}\int_{0}^{1}u^{a-1}v^{b-1}(1-u)^{c-1}(1-v)^{-e}(1-su) ^{b-e}
\left( 1-\dfrac{1}{t}v+\dfrac{s}{t}uv\right) ^{d-1}\,du\,dv
\label{4.9.21}
\end{multline}
Formula \eqref{4.9.21} is an integral representation of the $H_2$ function \cite[(3.4)]{7}. We get
\begin{equation}
I_{3}(s,t)=B(a,c)B(b,1-e)
s^{a+c-1}t^{d-1}H_2\left( e-b,a,1-d,b;a+c;s,-\dfrac{1}{t}\right)
\label{4.9.22}
\end{equation}
\subsection{The integral $I_4(s,t)$}
For $I_{4}(s,t)$ we look at $I_{2}(s,t)$. Interchanging in $I_{2}(s,t)$ the parameters $a$ and $b$, $c$ and $d$ and the variables $s$ and $t$ we get $I_{4}(s,t)$. So from \eqref{4.9.22} we get directly
\[
I_{4}(s,t)=B(a,1-e)B(b,d)
t^{b+d-1}s^{c-1}H_2\left( e-a,b,1-c,a;b+d;t,-\dfrac{1}{s}\right)
\]
\subsection{The integral $I_5(s,t)$}
For $I_{5}(s,t)$ we write the inner integral of \eqref{4.9.17e} as
\[
\int_{0}^{s}u^{a-1}(s-u)^{c-1}(1-v-u)^{-e}du=B(a,c)s^{a+c-1}(1-v) ^{-e}
\hyp21{a,e}{a+c}{\dfrac{s}{1-v}}
\]
with convergence conditions $\Re(a),\Re(c),\Re(c-e)>0$ and $\dfrac{s}{1-v}\leq 1$. Substitution in \eqref{4.9.17e} and replacing $v$ by $(1-s)v$ gives
\begin{multline}
	I_{5}(s,t) =B(a,c)s^{a+c-1}t^{d-1}(1-s)^{b} \\
	\times\int_{0}^{1}v^{b-1}\big( 1-(1-s)v\big)^{-e} 
	\left( 1-\dfrac{1-s}{t}v\right) ^{d-1}
	\hyp21{a,e}{a+c}{\dfrac{s}{1-(1-s)v}}dv
	\label{4.9.27}
\end{multline}

Because the computation of the integral is very long we only give the result from \cite[section 6.5.5.]{1}

\begin{align}
I_{5}(s,t)&
=B(c,e-c)B(b,c-e+1)s^{a-1}t^{d-1}(1-s)^{b+c-e} \nonumber \\
& \qquad \qquad \qquad \qquad \qquad \qquad  \times 
F_3\left( 
\begin{array}{c}
1-a,b;c,1-d \\ 
b+c-e+1%
\end{array}%
;\dfrac{s-1}{s},\dfrac{1-s}{t}\right) + \nonumber \\
&+B(a,c-e)B(b,d)s^{a+c-1}t^{b+d-1}F_P(e,b,a,b+d,a+c;t,s) - \nonumber \\
&-\dfrac{B(a,c-e)}{d}s^{a+c-e-1}t^{b-1}(s+t-1)^{d} \nonumber \\
& \qquad  \times
\widetilde{F_3}\left( 
\begin{array}{c}
1,1-b;e-a-c+1,d \\ 
d+1%
\end{array}%
;%
\begin{array}{c}
e \\ 
e-c+1%
\end{array}%
;\dfrac{s+t-1}{s},\dfrac{s+t-1}{t}\right)
\label{4.9.57}
\end{align}
with convergence conditions for the $F_P$ function: $\left\vert t\right\vert <1$, $\left\vert s-1\right\vert <1$. The convergence conditions of the $F_3$ function are: $\left\vert \dfrac{s-1}{s}\right\vert <1$ and $\left\vert \dfrac{1-s}{t}\right\vert<1$. This gives a region of convergence $\dfrac{1}{2}<s<1$. But after analytical continuation (see \eqref{4.P.20}) this function is well-defined for $\dfrac{s-1}{s}<1$ and $\dfrac{1-s}{t}<1$. For the extended $F_3$ function we have with analytical continuation convergence for the same region as for the $F_3$ function (see appendix B.2).This includes the $(s,t)$ domain on which the integral $I_5(s,t)$ is defined. For the parameters we have the convergence conditions: $0<\Re(a),\Re(b),\Re(c),\Re(d),\Re(e)<\Re(c)+1$.
 
\subsection{The integral $I_6(s,t)$}
For $I_6(s,t)$ we write the inner integral of \eqref{4.9.17f} as
\[
\int_{0}^{1-v}u^{a-1}(s-u)^{c-1}(1-v-u)^{-e}du= B(a,1-e)s^{c-1}(1-v)^{a-e}\hyp21{a,1-c}{a-e+1}{\dfrac{1-v}{s}}
\]
with convergence conditions $\Re(a),\Re(1-e)>0$ and $\dfrac{1-v}{s}<1$. 

Because the computation of the integral is also very long we only give the result from \cite[section 6.5.6.]{1}
\begin{align}
I_{6}(s,t)
&=\dfrac{B(a,c-e)}{d}s^{a+c-e-1}t^{b-1}(s+t-1)^{d} \nonumber \\
&\qquad\qquad \times 
\widetilde{ F_3}\left( 
\begin{array}{c}
1,1-b;e-a-c+1,d \\ 
d+1%
\end{array}%
;%
\begin{array}{c}
e \\ 
e-c+1%
\end{array}%
;\dfrac{s+t-1}{s},\dfrac{s+t-1}{t}\right) +\nonumber \\
&+B(1-e,e-c)B(d,1-e+c)s^{a-1}t^{b-1}(s+t-1)^{c+d-e} \nonumber \\
&\qquad\qquad \times F_3\left( 
\begin{array}{c}
1-a,1-b;c,d \\ 
c+d-e+1%
\end{array}%
;\dfrac{s+t-1}{s},\dfrac{s+t-1}{t}\right)
\label{4.9.58}
\end{align}
with convergence conditions:  $\Re(a),\Re(d),\Re(1-e),\Re(c-e+1)>0$. 

\subsection{The summation of the integrals $I_5(s,t)$ and $I_6(s,t)$}
When adding the functions $I_{5}(s,t)$ \eqref{4.9.57} and $I_{6}(s,t)$ \eqref{4.9.58}
the extended $F_3$ functions will be eliminated. We get
\begin{align}
I_{5}(s,t)+I_{6}(s,t)
&=B(a,c-e)B(b,d)s^{a+c-1}t^{b+d-1}F_P(e,b,a,b+d,a+c;t,s)+ \nonumber \\
&+B(c,e-c)B(b,c-e+1)s^{a-1}(1-s)^{b+c-e}t^{d-1} \nonumber \\
&\qquad \qquad \qquad \qquad \qquad F_3\left( 
\begin{array}{c}
1-a,b;c,1-d \\ 
b+c-e+1%
\end{array}%
;\dfrac{s-1}{s},\dfrac{1-s}{t}\right) + \nonumber \\
&+B(1-e,e-c)B(d,1-e+c)s^{a-1}t^{b-1}(s+t-1) ^{c+d-e} \nonumber \\
&\qquad \qquad \qquad \qquad F_3\left( 
\begin{array}{c}
1-a,1-b;c,d \\ 
c+d-e+1%
\end{array}%
;\dfrac{s+t-1}{s},\dfrac{s+t-1}{t}\right)
\label{4.9.59}
\end{align}
The power series of the functions $F_P$ \eqref{4.B.5} and $F_3$ \eqref{4.P.20} are convergent in the region $0<s<1,1-s<t<1$. For the first $F_3$ function we use the analytical continuation over the whole $(s,t)$ region. For the parameters we get the convergence conditions: $0<~\Re(a),\Re(b),\Re(c),\Re(d)$ and $\Re(e)<1$.
From \cite[p 441]{11} we get
\begin{align}
F_3\left( 
\begin{array}{c}
a_0,b_1;b_2,c_1 \\ 
c_2%
\end{array}%
;x_{1},x_{2}\right)
&=\dfrac{\Gamma (c_2) \Gamma(c_2-b_1-c_1) }{\Gamma( c_2-b_1)
	\Gamma(c_2-c_1) }(x_{1}) ^{-a_0}(x_{2})^{-b_1} \nonumber \\
&F_Q\left(b_2+c_1-c_2+1,b_2,c_1;1-a_0+b_2,c_1-b_1+1;\dfrac{1}{x_{1}},%
\dfrac{1}{x_{2}}\right) + \nonumber \\
&+\dfrac{\Gamma (c_2) \Gamma(b_1+c_1-c_2) }{\Gamma(b_1) \Gamma (c_1) }(x_{2})^{1-c_2}
(1-x_{2}) ^{c_2-b_1-c_1} \nonumber \\
&\qquad\quad F_3\left( 
\begin{array}{c}
a_0,1-b_1;b_2,1-c_1 \\ 
1-b_1-c_1+c_2%
\end{array}%
;\dfrac{x_{1}(x_{2}-1)}{x_{2}},1-x_{2}\right)
\label{4.9.60}
\end{align}
After application of \eqref{4.9.60} to the second $F_3$ function of \eqref{4.9.59} we get
\begin{align*}
I_{5}(s,t)+I_{6}(s,t)
&=B(a,c-e)B(b,d)s^{a+c-1}t^{b+d-1}F_P( e,b,a,b+d,a+c;t,s) +  \\
&+B(1-e,e-c)B(d,c-e+b)(s+t-1)^{a+b+c+d-e-2}  \\
&\qquad \qquad \qquad \qquad \qquad \qquad
F_Q\left( e,c,d;a+c,b+d;\dfrac{s}{s+t-1},\dfrac{t}{s+t-1}\right) + \\
&+\Gamma (e-c)\Gamma (c-e+1)\left[\dfrac{\Gamma (b)\Gamma (c) }{\Gamma(e)\Gamma(c+b-e+1) }+
\dfrac{\Gamma(1-e)\Gamma(e-b-c) }{\Gamma(1-c)\Gamma(1-b) }\right]  \\
&\qquad \qquad \qquad \qquad s^{a-1}(1-s) ^{b+c-e}t^{d-1}F_3\left( 
\begin{array}{c}
1-a,b;c,1-d \\ 
b+c-e+1%
\end{array}%
;\dfrac{s-1}{s},\dfrac{1-s}{t}\right)
\end{align*}
The function $F_Q$ is defined in \eqref{4.B.5a}. After some simplification of the terms within the square brackets we get
\begin{align*}
I_{5}(s,t)+I_{6}(s,t)
&=B(a,c-e)B(b,d)s^{a+c-1}t^{b+d-1}F_P( e,b,a,b+d,a+c;t,s) +  \\
&+B(1-e,e-c)B(d,c-e+b)s^{a+c-1}t^{b+d-1}(s+t-1) ^{-e}  \\
&\qquad \qquad \qquad \qquad \qquad \qquad  
F_Q\left( e,c,d;a+c,b+d;\dfrac{s}{s+t-1},\dfrac{t}{s+t-1}\right) +  \\
&+B(e-b-c,c)B(1-e,b)s^{a-1}\left( 1-s\right) ^{b+c-e}t^{d-1}  \\
&\qquad \qquad \qquad \qquad \qquad \qquad \qquad \qquad F_3\left( 
\begin{array}{c}
1-a,b;c,1-d \\ 
b+c-e+1%
\end{array}%
;\dfrac{s-1}{s},\dfrac{1-s}{t}\right) 
\end{align*}
The power series of the $F_P$ and the $F_Q$ functions are convergent in the region $\{0<s<1 \wedge 1-s<t<1\}$. For the $F_3$ function we use an analytical continuation so we get the convergence in the right region. For the parameters we get the convergence conditions: $0<\Re(a,b,c,d,e)<1$.

\subsection{Summary of the integrals $I_{1..6}(s,t)$}
We give a summary of the functions $I_{1..9}(s,t)$
\begin{equation}
I_{1}(s,t)=B(a,c)B(b,d)s^{a+c-1}t^{b+d-1}F_2(e;a,b;a+c,b+d;s,t)
\label{4.I78.a}
\end{equation}
with $0<s<1,0<t<1-s$ and $0<\Re(a)<\Re(a+c)$, $0<\Re(b)<\Re(b+d)$.
\[
I_{2}(s,t)=B(a,b,1-e)s^{c-1}t^{d-1}
F_3\left(a,b;1-c,1-d;a+b+1-e;\dfrac{1}{s},\dfrac{1}{t}\right)
\] 
with $1<s<\infty ,1<t<\infty$ and $\Re(a),\Re(b),\Re(1-e)>0$.
\begin{equation}
I_{3}(s,t)=B(a,c)B(b,1-e)s^{a+c-1}t^{d-1}
H_2\left( e-b,a,1-d,b;a+c;s,-\dfrac{1}{t}\right)
\label{4.I78.c}
\end{equation}
with $0<s<1<t<\infty$ and $\Re(a),\Re(b)>0,\Re(c)>0$.
\begin{equation}
I_{4}(s,t)=B(a,1-e)B(b,d)t^{b+d-1}s^{c-1}
H_2\left( e-a,b,1-c,a;b+d;t,-\dfrac{1}{s}\right)
\label{4.I78.d}
\end{equation}
with $0<t<1<s<\infty$ and $\Re(a),\Re(b),\Re(d)>0$.
\begin{align}
I_{5}(s,t)+I_{6}(s,t)
&=B(a,c-e)B(b,d)s^{a+c-1}t^{b+d-1}F_P(e,b,a,b+d,a+c;t,s) + \nonumber \\
&+B(c,e-c)B(b,c-e+1)s^{a-1}(1-s) ^{b+c-e}t^{d-1} \nonumber \\
&\qquad \qquad \qquad \qquad \qquad 
F_3\left( 
\begin{array}{c}
1-a,b;c,1-d \\ 
c+b-e+1%
\end{array}%
;\dfrac{s-1}{s},\dfrac{1-s}{t}\right) + \nonumber \\
&+B(1-e,e-c)B(d,1-e+c)t^{b-1}s^{a-1}(s+t-1) ^{c+d-e} \nonumber \\
&\qquad \qquad \qquad \qquad 
F_3\left( 
\begin{array}{c}
1-a,1-b;c,d \\ 
c+d-e+1%
\end{array}%
;\dfrac{s+t-1}{s},\dfrac{s+t-1}{t}\right)
\label{4.I78.e}
\end{align}
with $0<s<1,\ 1-s<t<1$ and $0<\Re(a),\Re(c),\Re(d),\Re(e)<1$.
\\[3mm]
The functions $F_2$, $F_3$, $H_2$ and $F_P$ are all five-parametric double-hypergeometric functions which are solutions of a system of two partial differential equations of second order  \cite{9}. This system is associated with the Appell function $F_2(a,b_1,b_2,c_1,c_2;x,y)$ \cite[5.9(10)]{28}. Note that the orthogonal polynomials 
\eqref{4.2.16} and \eqref{4.2.17} used for the fractional orthogonal derivative are also expressed in terms of $F_2$ functions. 	
\subsection{Remarks}
\begin{remark}\
	It is well-known \cite[5.9(10)]{28} that the $F_2(a_0;b_1,b_2;c_1,c_2;x,y)$ function is a solution of the following system of partial differential equations
	\begin{equation} 
	x(1-x) r-xys+\Big[ c_1-(a_0+b_1+1) x \Big] p-b_1yq-a_0b_1F_2(x,y) =0
	\label{4.9.75a}
	\end{equation}
	\begin{equation}
	y(1-y) t-xys+\Big[ c_2-(a_0+b_2+1) y \Big] q-b_2xp-a_0b_2F_2(x,y) =0
	\label{4.9.75b}
	\end{equation}
	with
	\[
	p=\dfrac{\partial F_2(x,y) }{\partial x} \ \ \
	q=\dfrac{\partial F_2(x,y)}{\partial y} \ \ \
	r=\dfrac{\partial	^{2}F_2(x,y) }{\partial x^{2}} \ \ \
	s=\dfrac{\partial^{2}F_2(x,y)}{\partial x \partial y} \ \ \
	t=\dfrac{\partial^{2}F_2(x,y)}{\partial y^{2}}
	\]
	Olsson treats all 36 explicit solutions of \eqref{4.9.75a} and \eqref{4.9.75b} \cite{9}. Using the substitutions $a_0=e-a-b-c-d+2,\ b_1=1-c,\ b_2=1-d,\ c_1=2-a-c,\ c_2=2-b-d$ and $x\rightarrow s,\ y\rightarrow t$ he gives among others the following solutions of \eqref{4.9.75a} and \eqref{4.9.75b}
	\begin{equation}
	s^{a+c-1}t^{b+d-1}F_2\left( 
	\begin{array}{c}
	e;a,b \\ a+c,b+d%
	\end{array}%
	;s,t\right)
	\label{4.9.78a}
	\end{equation}
	\begin{equation}
	s^{a+c-1}F_2\left( 
	\begin{array}{c}
	e-b-d+1;a,1-d \\ 
	a+c,2-b-d%
	\end{array}%
	;s,t\right)
	\label{4.9.78b}
	\end{equation}
	\begin{equation}
	s^{c-1}t^{d-1}F_3\left( 
	\begin{array}{c}
	a,b;1-c,1-d \\ 
	a+b+1-e%
	\end{array}%
	;\dfrac{1}{s},\dfrac{1}{t}\right)
	\label{4.9.78c}
	\end{equation}
	\begin{equation}
	t^{b+d-1}F_P(e-c-a+1,b,1-c,b+d,2-a-c;t,s)
	\label{4.9.78d}
	\end{equation}
	\begin{equation}
	s^{a-1}t^{d-1}(1-s) ^{b+c-e}F_3\left( 
	\begin{array}{c}
	1-a,b;c,1-d \\ 
	b+c-e+1%
	\end{array}%
	;\dfrac{s-1}{s},\dfrac{1-s}{t}\right)
	\label{4.9.78e}
	\end{equation}
	\begin{equation}
	s^{a-1}t^{b-1}(1-s-t) ^{c+d-e}F_3\left( 
	\begin{array}{c}
	1-a,1-b;c,d \\ 
	c+d-e+1%
	\end{array}%
	;\dfrac{s+t-1}{s},\dfrac{s+t-1}{t}\right)
	\label{4.9.78f}
	\end{equation}
	\begin{equation}
	s^{a+c-1}t^{d-1}
	H_2\left( e-b,a,1-d,b;a+c;s,-\dfrac{1}{t}\right)
	\label{4.9.78g}
	\end{equation}
	
	\ 
	
	For the integral $I_{1}(s,t)$ we see that the function $s^{a+c-1}t^{b+d-1}F_2\left(
	\begin{array}{c}
	e;a,b \\ 
	a+c,b+d%
	\end{array}%
	;s,t\right) $ equals the solution \eqref{4.9.78a}.
	
	\ 
	
	For the integral $I_{2}(s,t)$ we see that the function $s^{c-1}t^{d-1}F_3\left( a,b;1-c,1-d;a+b+1-e;\dfrac{1}{s},\dfrac{1}{t}\right)$  equals the solution \eqref{4.9.78c}.
	
	\
	
	For the integral $I_{3}(s,t)$ we see that the function $s^{a+c-1}t^{d-1}H_2\left( e-b,a,1-d,b;a+c;s,-\dfrac{1}{t}\right)$  equals the solution \eqref{4.9.78g}.
	
	\ 
	
	To compute $I_{4}(s,t)$ we can interchange in $I_{3}(s,t)$ the parameters $a$ and $b$, $c$ and $d$ and the variables $s$ and $t$. Then we see that \eqref{4.I78.d} is a solution of \eqref{4.9.75a} and \eqref{4.9.75b}.
	
	\ 
	
	For the integral $I_{5}(s,t)+I_{6}(s,t)$ we have the functions
	\begin{equation}
	L_{1}(s,t) =s^{a+c-1}t^{b+d-1}F_P(e,b,a,b+d,a+c;t,s)
	\label{4.9.81a}
	\end{equation}
	
	\[
	L_{2}(s,t) =s^{a-1}(1-s)^{b+c-e}t^{d-1}F_3\left( 
	\begin{array}{c}
	1-a,b;c,1-d \\ 
	b+c-e+1%
	\end{array}%
	;\dfrac{s-1}{s},\dfrac{1-s}{t}\right)
	\]
	
	\[
	L_{3}\left(s,t\right) =s^{a-1}t^{b-1}(s+t-1) ^{c+d-e}F_3\left( 
	\begin{array}{c}
	1-a,1-b;c,d \\ 
	c+d-e+1%
	\end{array}%
	;\dfrac{s+t-1}{s},\dfrac{s+t-1}{t}\right)
	\]
	For the function $L_{1}(s,t)$ we use the solution \cite[(26)]{9} and interchange $b_1\longleftrightarrow b_2,c_1\longleftrightarrow c_2$ and $x_{1}\longleftrightarrow
	x_{2}$. Then the solution becomes
	\[
	\left( x_{2}\right) ^{1-c_2}F_P(a_0-c_2+1,b_1,b_2c_2+1,c_1,2-c_2;x_{2},x_{1})
	\]
	Application to \eqref{4.9.81a} gives as result
	\[
	t^{b+d-1}F_P(e-c-a+1,b,1-c,b+d,2-a-c;t,s)
	\]
	This function is equal \eqref{4.9.78d}. So $L_{1}(s,t)$ is a solution of \eqref{4.9.78a} and \eqref{4.9.78b}.
	
	\ 
	
	The function $L_{2}(s,t)$ equals the solution \eqref{4.9.78e}.
	
	\ 
	
	The function $L_{3}(s,t)$ equals the solution \eqref{4.9.78f}.
\end{remark}

\begin{remark}
	It is clear that the inner double integrals in \eqref{4.9.11} are in fact double integrals transformations which transform the solutions of the $F_2$ equations \eqref{4.9.75a} and \eqref{4.9.75b} into a solution of the $F_2$ equations with other parameters.
\end{remark}

\begin{remark}\
	When looking at Figure \ref{Figure 6.1} case Va we can consider the integral when the arguments are on the boundaries. Then we distinguish four possibilities.
	
	\underline{I.\qquad $s=1$.}
	
	For $s=1$ it can be shown that $I_{5}(1,t)+I_{6}(1,t)=I_{4}(1,t)$. 
	\\[4mm]	
	\underline{II.\qquad $t=1$.}
	
	For $t=1$ it can be shown that $I_{5}(s,1)+I_{6}(s,1)=I_{3}(s,1)$. 
	\\[4mm]		
	\underline{III.\qquad $s=t=1$.}
	
	From \eqref{4.9.27} it follows $I_{5}(1,t)=0$ and $I_{6}(1,t)=I_{4}(1,t)$. For $s=t=1$ it can be shown that $I_{6}( 1,1)=I_{2}(1,1)=I_{3}(1,1)$. 
	\\[4mm]	
	\underline{IV.\qquad $s+t=1$.}
	
	For $s+t=1$ it can be shown that $I_{5}(s,1-s)+I_{6}(s,1-s)=I_{1}(s,1-s)$.
\end{remark}

\begin{remark}\
	When looking at case V in Figure \ref{Figure 6.1} we see directly that when interchanging $a\leftrightarrow b$, $c\leftrightarrow d$, $s\leftrightarrow t$ and the integration variables $u\leftrightarrow v$ we get
	\begin{equation}
	I_{5}(a,b,c,d,e;s,t)+I_{6}(a,b,c,d,e;s,t) =I_{5}(b,a,d,c;e;t,s) +I_{6}(b,a,d,c;e;t,s)
	\label{4.9.124}
	\end{equation}
	From the symmetry there follows that $I_5(s,t)+I_6(s,t)$ can be written as half the sum of the two sides of \eqref{4.9.124}.
	
	Looking at equation \eqref{4.I78.e} we see that only the second $F_3$ function in the right-hand side is symmetrical in the parameters and the arguments. Now write \eqref{4.9.124} as
	\begin{align*}
	I_{5}(a,b,c,d,e;s,t)+I_{6}(a,b,c,d,e;s,t)
	&=\lambda\big( I_{5}(a,b,c,d,e;s,t)+I_{6}(a,b,c,d,e;s,t) \big) +\\
	&+(1-\lambda)\big( I_{5}(b,a,d,c;e;t,s) +I_{6}(b,a,d,c;e;t,s) \big)
	\end{align*}
	Now we can choose $\lambda$ so that the second $F_3$ functions will be eliminated. There results
	\begin{align*}
	I_{5}(s,t)+I_{6}(s,t)
	&=A_{1}\;s^{a+c-1}t^{b+d-1}F_P( e,b,a,b+d,a+c;t,s)+ \\
	&+A_{2}\;s^{a+c-1}t^{b+d-1}F_P( e,a,b,a+c,b+d;s,t)+ \\
	&+A_{3}\;t^{b-1}(1-t)^{a+d-e}s^{c-1} 
	F_3\left( 
	\begin{array}{c}
	d,a;1-b,1-c \\ 
	a+d-e+1%
	\end{array};\dfrac{t-1}{t},\dfrac{1-t}{s}\right)+ \\
	&+A_{4}\;s^{a-1}(1-s)^{b+c-e}t^{d-1}
	F_3\left( 
	\begin{array}{c}
	c,b;1-a,1-d \\ 
	b+c-e+1%
	\end{array};\dfrac{s-1}{s},\dfrac{1-s}{t}\right)
	\end{align*}
	with
	\[
	A_1=\dfrac{\Gamma(a)\Gamma(b)\Gamma(e)\Gamma(c-d)\Gamma(d-c+1)\Gamma(1-e)}
	{\Gamma(a+c-e)\Gamma( b+d) \Gamma(1-d)\Gamma(e-c+1)}
	\]
	\[
	A_2=\dfrac{\Gamma(a)\Gamma(b)\Gamma(e)\Gamma(d-c)\Gamma(c-d+1)\Gamma(1-e)}
	{\Gamma( b+d-e)\Gamma(a+c)\Gamma(1-c)\Gamma(e-d+1)}
	\]
	\[
	A_3=\dfrac{\Gamma(a)\Gamma(d)\Gamma(c-d)\Gamma(d-c+1)\Gamma(1-e)}
	{\Gamma(c)\Gamma(1-c)\Gamma(a+d-e+1)} \qquad
	\]
	\[
	A_4=\dfrac{\Gamma(b)\Gamma(c)\Gamma(d-c)\Gamma(c-d+1)\Gamma(1-e)}
	{\Gamma(d)\Gamma(1-d)\Gamma(c+b-e+1)} \qquad
	\]
	Application of Olsson \cite[(26)]{9} gives
	\begin{align*}
	I_{5}(s,t)+I_{6}(s,t)
	&=A_{1}\;s^{a+c-1}t^{b+d-1}(1-t)^{-e}F_P\left( e,d,a,b+d,a+c;\dfrac{t}{t-1},\dfrac{s}{1-t}\right)+ \\
	&+A_{3}\;s^{c-1}t^{b-1}(1-t)^{a+d-e} 
	F_3\left( 
	\begin{array}{c}
	d,a;1-b,1-c \\ 
	a+d-e+1%
	\end{array};\dfrac{t-1}{t},\dfrac{1-t}{s}\right)+ \\
	&+A_{2}\;s^{a+c-1}t^{b+d-1}(1-s)^{-e}F_P\left( e,c,b,a+c,b+d;\dfrac{s}{s-1},\dfrac{t}{1-s}\right)+ \\
	&+A_{4}\;s^{a-1}t^{d-1}(1-s)^{b+c-e}
	F_3\left( 
	\begin{array}{c}
	c,b;1-a,1-d \\ 
	b+c-e+1%
	\end{array};\dfrac{s-1}{s},\dfrac{1-s}{t}\right)
	\end{align*}
	Note that the arguments of the $F_P$ functions are the inverses of the arguments of the $F_3$ functions.
\end{remark}
\

\

\textbf{\Large Appendices:  Hypergeometric functions of two variables}

\

When deriving formulas for the two-dimensional fractional orthogonal derivatives we encounter various hypergeometric functions of two variables, some of which are beyond the familiar Appell and Horn hypergeometric functions. Here we give an overview of the functions we need in this thesis.

\begin{appendices}
	\section{The Pochhammer symbol}
	
We have already frequently used the Pochhammer symbol, which we defined for
$a\in\mathbb{C}$ by
\[
(a)_i:=a(a+1)\ldots(a+i-1)=\frac{\Gamma(a+i)}{\Gamma(a)}\quad(i=0,1,\ldots),
\]
and more generally by
\[
(a)_b:=\frac{\Gamma(a+b)}{\Gamma(a)}\quad(b\in\mathbb{C}).
\]
We will use this second definition in particular for $b\in\mathbb{Z}$. Then
\[
(a)_{-i}=\frac{(-1)^i}{(1-a)_i}\quad(i\in\mathbb{Z}).
\]
Further properties for $i,j\in\mathbb{Z}$ are:
\[
\Gamma(a+i)=\Gamma(a)(a) _{i}
\]
\[
\Gamma(a-i) =(-1) ^{i}\dfrac{\Gamma(a)}{(1-a) _{i}}
\]
\[
(a+i)_j=\dfrac{(a)_{i+j}}{(a)_i}=\dfrac{(a+j)_i (a)_j}{(a)_i}
\]
\[
(a-i)_{j}=(-1)^{i}(a)_{j-i}(1-a)_i=(-1)^{j}\dfrac{(1-a)_{i}}{(1-a)_{i-j}}
=(a)_j\dfrac{(1-a)_i}{(1-a-j)_i}
\]
There are a lot more. See for example \cite[Appendix I]{17}.
	
\section{Definitions}
	
For the basics of hypergeometric functions of two variables see for instance
\cite[Ch. 9]{15}, \cite[sections 5.7-5.12]{28}, \cite[Ch. 1]{16}, \cite[Ch. 8]{17}.
	
Appell \cite{5}, \cite{8} defines four functions $F_1, F_2, F_3, F_4$ in two variables which are generalizations of the \,$_{2}F_1$\, hypergeometric function. We use only the first three functions \cite[5.7]{28}.
	\[
	F_1(a;b_1,b_2;c;x,y) =\sum\limits_{i=0}^{\infty}\sum\limits_{j=0}^{\infty }
	\dfrac{(a) _{i+j}(b_1)_{i}(b_2)_{j}}{(c) _{i+j}}\dfrac{1}{i!\,j!}x^{i}y^{j}
	\]
	with convergence region: 
	$\left\{\left\vert x\right\vert <1\wedge \left\vert y\right\vert <1\right\}$.
	\[
	F_2(a;b_1,b_2;c_1,c_2;x,y)=\sum\limits_{i=0}^{\infty }\sum\limits_{j=0}^{\infty }
	\dfrac{(a)_{i+j}(b_1)_{i}(b_2)_{j}}{(c_1) _{i}(c_2)_{j}}\dfrac{1}{i!\,j!}x^{i}y^{j}
	\]
	This power series has as convergence region: $\{\vert x \vert +\vert y \vert <1\}$. From the known transformations of the $F_2$ function \cite[5.11. (6),(7),(8)]{28} and from the union of the four partially overlapping regions of convergence (or from the integral representation of the $F_2$ function \cite[5.8(2)]{28}) we get a region of unique analytic convergence $\left\{ x<1\wedge y<1\wedge x+y<1\right\}$. 
	\begin{equation}
	F_3(a_{1},a_{2};b_1,b_2;c;x,y)=\sum\limits_{i=0}^{\infty}\sum\limits_{j=0}^{\infty }
	\dfrac{(a_{1})_{i}(a_{2})_{j}(b_1)_{i}(b_2) _{j}}{(c) _{i+j}}\dfrac{1}{i!\,j!}x^{i}y^{j}
	\label{4.P.20}
	\end{equation}
	with convergence region: $\left\{\vert x\vert <1\wedge \vert y\vert <1\right\}$. With the integral representation of the $F_3$ function \cite[5.8(3)]{28}\footnote{Note that the power of the factor $(1-u-v)$ should be $\gamma-\beta-\beta'.$} we get a region of unique analytical continuation of $\{x<1,\,y<1\}$.
	
	Another integral representation is
	\begin{multline}
	F_3(a_{1},a_{2};b_1,b_2;c;x,y)=
	\dfrac{\Gamma(c)}{\Gamma(a_1)\Gamma(a_2)\Gamma(c-a_1-a_2)} \\
	\times\int_0^1\int_0^1 u^{a_1-1}v^{a_2-1}(1-u)^{c-a_2-a_1-1}(1-v)^{c-a_2-1}
	(1-y\,v)^{-b_2}(1-x\,u+x\,u\,v)^{-b_1}dudv
	\label{4.P.20a}
	\end{multline}
	with $\Re(a_1,a_2,c-a_2-a_1,c-a_2)>0$. 
	
	We define the {\em extended} $F_3$ function as
	\begin{equation}
	\widetilde{F_3}\left(\begin{array}{c}a_{1},a_{2},b_1,b_2 \\ 
	c\end{array};\begin{array}{c}d_{1} \\ 
	d_{2}\end{array};x,y\right)
	=\sum\limits_{i=0}^{\infty}\sum\limits_{j=0}^{\infty }\dfrac{(a_{1})_{i}(a_{2})_{j}(b_1) _{i}(b_2)_{j}}{(c)_{i+j}}\dfrac{(d_1)_i}{(d_2)_i}\dfrac{1}{i!\,j!}x^{i}y^{j}
	\label{4.P.21}
	\end{equation}
	This extended $F_3$ function is of order three (see below) and can be written as a triple integral
	\begin{multline*}
	\widetilde{F_3}\left(\begin{array}{c}a_{1},a_{2},b_1,b_2 \\ 
	c\end{array};\begin{array}{c}d_{1} \\ 
	d_{2}\end{array};x,y\right)
	=\dfrac{\Gamma(c)}{\Gamma(a_1)\Gamma(a_2)\Gamma(c-a_1-a_2)}
	\dfrac{\Gamma(d_2)}{\Gamma(d_1)\Gamma(d_2-d_1)} \\
	\qquad\qquad\times\int_0^1\int_0^1\int_0^1 u^{a_1-1}v^{a_2-1}w^{d_1-1}
	(1-u)^{c-a_1-1}(1-v)^{c-a_1-a_2-1}(1-w)^{d_2-d_1-1} \\
	\times(1-yv+yuv)^{-b_2}(1-x\,u\,w)^{-b_1}
	dudvdw
	\end{multline*}
	Then we have the convergence conditions for the parameters: $\Re(c)>\Re(a_1)>0,\Re(c-~a_1)>\Re(a_2)>0,\Re(d_2)>\Re(d_1)$. For the arguments we have the analytical continuation of convergence: $\{x<1 \wedge y<1\}$.
	
	Horn \cite{12} more generally calls a double power series
	\[
	\sum\limits_{i,j=0}^{\infty }A\left( i,j\right) x^{i}y^{j}
	\]
	{\em hypergeometric} if the two quotients $\dfrac{A(i+1,j) }{A(i,j) }$ and $ \dfrac{A(i,j+1) }
	{A(i,j) }$ are rational functions of $i$ and $j$. Write the two rational functions as quotients of polynomials in $i$ and $j$ without common factors.
	
	\
	
	\qquad $\dfrac{A(i+1,j)}{A(i,j)}=\dfrac{F(i,j) }{F^{\prime }(i,j) }\qquad $and\qquad
	$\dfrac{A(i,j+1)}{A(i,j) }=\dfrac{G(i,j)}{G^{\prime	}(i,j)}$
	
	\
	
	In addition it is assumed that $F'(i,j)$ contains the factor $i+1$ and $G'(i,j)$ the 
	factor $j+1$. Then the highest degree in $i,j$ of the four polynomials $F,F^{\prime},G,G^{\prime }$ is called the $order$ of the hypergeometric series. Horn also gives a rule (see \cite[Section 5.7.2]{28}) how to determine the region of convergence of the power series from the two above quotients. Then Horn classifies the hypergeometric functions of order two. He gives a list of 34 such functions. The first four items in the list are the Appell functions $F_1,\ F_2,\ F_3$ and $F_{4}$.
	
	\
	
	Number nine in the list of Horn is the $H_2$ function.
	\begin{equation}
	H_2(a,b_1,b_2,c_1,c_2;x,y)=\sum\limits_{i=0}^{\infty}
	\sum\limits_{j=0}^{\infty }\dfrac{(a)_{i-j}(b_1)_{i}(b_2)_{j}(c_1)_{j}}{(c_2) _{i}}
	\dfrac{1}{i!\,j!}x^{i}y^{j}
	\label{4.P.1}
	\end{equation}
	with convergence condition (determined by Horn's rule):
	\begin{equation}
	\left\{|x|<1,|y|<(|x|+1)^{-1}\right\}
	\label{4.P.1a}
	\end{equation}
	In \cite[Section 3]{7} Diekema and Koornwinder showed that 
	\[
	\{ (x<0\wedge (x-1) y<1)\vee (0 \leq x<1\wedge y>-1)\}
	\]
	is a region in $\mathbb{R}^2$ where $H_2$ has unique analytical continuation.
	\begin{figure}[ht]
		\centering
		\parbox{5cm}{
			\includegraphics[width=5cm]{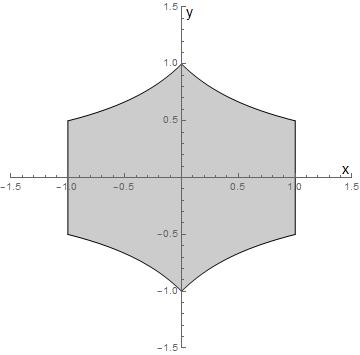}
			\caption{Convergence region of the $H_2$ function}
			\label{fig:1}}
		\qquad
		\begin{minipage}{5cm}
			\includegraphics[width=5cm]{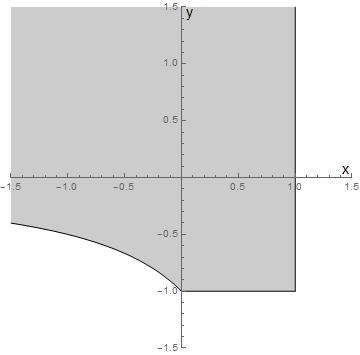}
			\caption{Region of analytical continuation of $H_2$}
			\label{fig:2}
		\end{minipage}
	\end{figure}
	
	In \cite[section 5.7]{28} for each function in Horn's list a system of two pde's is given which has this function as a solution. This follows a method already indicated by Horn in \cite{14}, \cite{12}.
	
	A comprehensive list of the solutions of the system of pde's \eqref{4.9.75a},\eqref{4.9.75b} for the $F_2$ function was given by Olsson, \cite[p.1289, table I]{9}. Apart from $F_2,\ F_3$ and $H_2$ there are solutions expressed in terms of functions $F_P,\ F_Q$ and $F_R$ which do not occur in Horn's list \cite{11}. Here $F_P$ and $F_Q$ have order three, while $F_R$ is not even of hypergeometric type. In his paper he also mentioned an $F_{PR}$ function \cite[(57)]{9}.
	
	\
	
	Olsson defines his $F_P$ function as:
	\begin{equation}
	F_P(a,b_1,b_2,c_1,c_2;x,y)=y^{-a}\sum\limits_{i=0}^{\infty }\sum\limits_{j=0}^{\infty }
	\dfrac{(a)_{i+j}(a-c_2+1)_{i+j}(b_1)_{i}}{(a+b_2-c_2+1)_{i+j}(c_1)_{i}}\dfrac{1}{i!\,j!}
	\left( \dfrac{x}{y}\right) ^{i}\left( \dfrac{y-1}{y}\right) ^{j}
	\label{4.P.11}
	\end{equation}
	with convergence region: $\{|xy^{-1}|+|1-y^{-1}|<1\}$. In \cite{11} (in a different notation) he gives
	\begin{equation}
	F_P(a,b_1,b_2,c_1,c_2;x,y)=
	y^{-a}\sum\limits_{i=0}^{\infty }\sum\limits_{j=0}^{\infty }
	\dfrac{(a)_{i+j}(a-c_2+1)_i(b_2)_j}{(a+b_2-c_2+1)_{i+j}}\dfrac{(b_1)_i}{(c_1)_i}\dfrac{1}{i!\,j!}
	x^{i}(1-y) ^{j}
	\label{4.B.5}
	\end{equation}
	with convergence region $\{|x|<1,|y-1|<1\}$. In \cite[Theorem 4.1]{7} we showed that
	$\{x<1\wedge y>0\}$ is a region in $\mathbb{R}^2$ where $F_P$ has unique analytical continuation.
	
	\
	
	Olsson \cite{11}, \cite{9} defined his $F_Q$ function as 
	\begin{multline}
	F_Q( a,b_1,b_2,c_1,c_2;x,y)=x^{-b_1}y^{b_1-a} \\
	\times\sum\limits_{i=0}^{\infty}\sum\limits_{j=0}^{\infty }\dfrac{(b_1+c_2-b_2-a)_{i-j}
		(b_1)_i(b_1-c_1+1)_i}{(b_1-a+1)_{i-j}(b_1+c_2-a)_{i-j}}
	\dfrac{1}{i!\,j!}\left( \dfrac{y}{x}\right)^i\left(\dfrac{1-y}{y}\right)^j
	\label{4.B.5a}
	\end{multline}
	where the convergence conditions (by Horn's rule) are the same as for the function 
	
	$H_2\left(1-\dfrac{1}{y},-\dfrac{y}{x}\right)$ \cite[(3)]{9}. So the convergence conditions following \eqref{4.P.1a} are:
	\[
	\left\{(x,y) \in \RR^{2}\text{ }|\text{ }|y-1|<|y|, |y-1|+|y|<|x|\right\}
	\]
	
	\begin{figure}[ht]
		\centering
		\includegraphics[width=5cm]{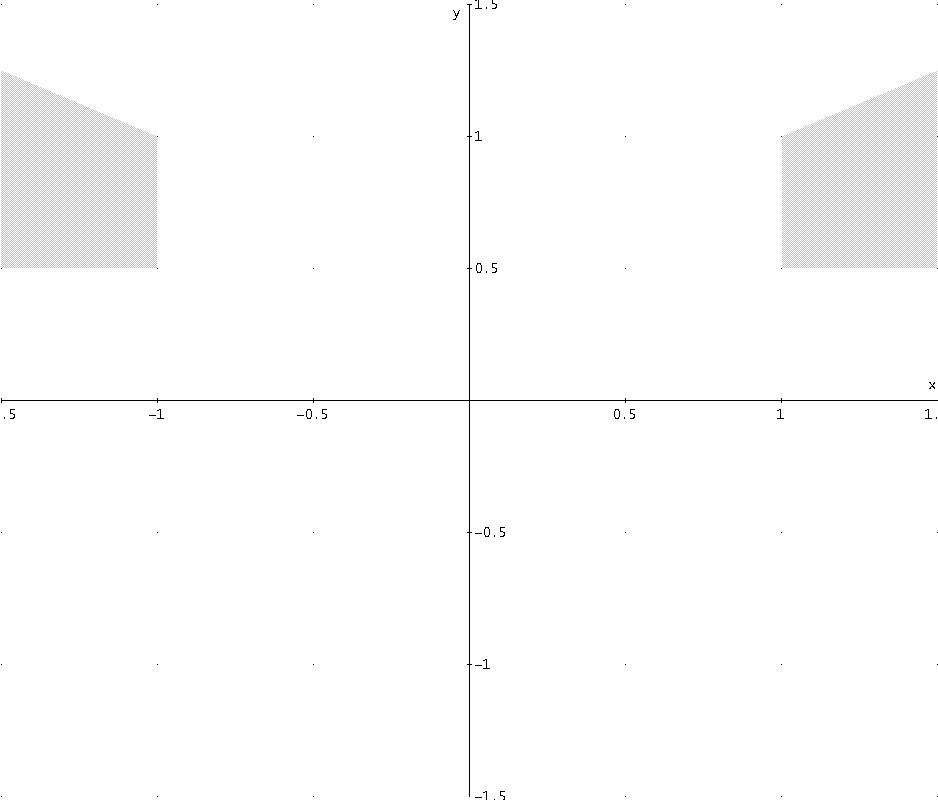}
		\caption{Convergence region of the $F_Q$ function}
		\label{fig:3}
	\end{figure}
	
	Olsson gives formula \cite[(56)]{9}
	\begin{multline}
	F_P(a,b_1,b_2;c_1,c_2;x,y) =
	F_{PR}(a,b_1,b_2;c_1,c_2;x,y)+ \\
	+\dfrac{\Gamma(a+b_2-c_2+1)\Gamma(a+b_1-c_1-b_2)\Gamma(c_1)}
	{\Gamma(a)\Gamma(b_1)\Gamma(a-c_2+1)}
	x^{b_1-c_1}y^{-b_2}(1-x)^{c_1-b_1+b_2-a} \\
	\times F_3\left( 
	\begin{array}{c}
	1-b_1,b_2;c_1-b_1,b_2-c_2+1 \\ 
	c_1-b_1+b_2-a+1%
	\end{array}%
	;\dfrac{x-1}{x},\dfrac{1-x}{y}\right)
	\label{4.B.6}
	\end{multline}
	where $F_{PR}$ is the part of the function $F_P$ which is regular at $(1,1)$. For the function $F_{PR}$ he gave the formula \cite[(57)]{9}
	\begin{multline}
	F_{PR}(a,b_1,b_2;c_1,c_2;x,y)=
	\dfrac{\Gamma(a+b_2-c_2+1)\Gamma(c_1-b_1+b_2-a)\Gamma(c_1)}
	{\Gamma(a)\Gamma(c_1-b_1+b_2-c_2+1)\Gamma(c_1+b_2-a)} \\
	\qquad\qquad\qquad\qquad \times\sum\limits_{i=0}^{\infty }\sum\limits_{j=0}^{\infty} 
	\dfrac{(a-c_2+1)_i(b_1)_i(b_2)_j}{(a+b_1-c_1-b_2+1)_i}\dfrac{1}{i!\,j!}(1-x)^i(1-y)^j\\
	\times\hyp32{b_2-c_2+1,c_1-b_1+b_2-a-i,c_1-a-j}{c_1-b_1+b_2-c_2+1,c_1+b_2-a}{1}
	\label{4.B.6a}
	\end{multline}
	where the $_3 F_2$ function is absolute convergent if $\Re(a)>0$. He then referred to the paper \cite[(12)]{93} in which the convergence region of the double summation is given as
	\[
	|x-1|+|y-1|<1
	\]
	The region of convergence of the $F_3$ function is:
	\[
	\left\{\left|\dfrac{x-1}{x}\right|<1 \wedge \left|\dfrac{1-x}{y} \right|<1\right\}
	\]
	The point $(1,1)$ is inside both regions. 
	\begin{figure}[ht]
		\centering
		\parbox{5cm}{
			\includegraphics[width=5cm]{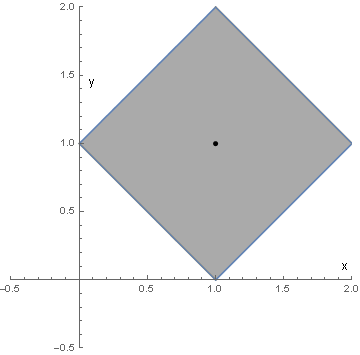}
			\caption{Convergence region of the $F_{PR}$ function}
			\label{fig:5}}
		\qquad
		\begin{minipage}{5cm}
			\includegraphics[width=6cm]{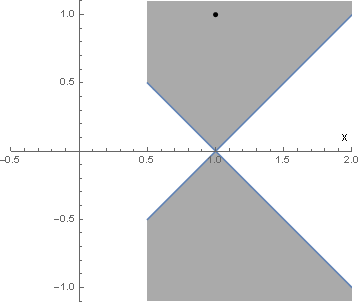}
			\caption{Convergence region of the $F_3$ function}
			\label{fig:4}
		\end{minipage}
	\end{figure}

\end{appendices}

\

\

\end{document}